\documentstyle[amscd,amssymb,12pt]{amsart}
\input xypic \xyoption{curve} 

\def\hcorrection#1{\advance\hoffset by #1 }
\def\vcorrection#1{\advance\voffset by #1 }

\vcorrection{-.7in}
\hcorrection{-.5in}

\newcommand{\B}[1]{{\bold#1}} 
\newcommand{\C}[1]{{\cal#1}} 

\theoremstyle{plain}
\newtheorem{th}{Theorem}[section]
\newtheorem{cor}{Corollary}[section]
\newtheorem{lem}{Lemma}[section]
\newtheorem{prop}{Proposition}[section]

\theoremstyle{definition}
\newtheorem{defin}{Definition}[section]

\theoremstyle{definition}
\newtheorem{example}{Example}[section]

\theoremstyle{remark}
\newtheorem{rem}{Remark}[section]

\numberwithin{equation}{section}

\begin{document}
\pagestyle{plain}
\addtolength{\footskip}{.3in}

\title{On ideals and homology in additive categories}
\author{Lucian M. Ionescu}
\address{Mathematics Department\\Kansas State University\\
             Manhattan, Kansas 66502}
\email{luciani@@math.ksu.edu}
\date{January, 2000}

\begin{abstract}
Ideals are used to define homological functors in additive categories.
In abelian categories the ideals corresponding to the usual
universal objects are principal, and the construction
reduces, in a choice dependent way, to homology groups.
Applications are considered: derived categories and functors.
\end{abstract}

\maketitle
\tableofcontents


\section{Introduction}\label{S:intro}
Categorification is by now a commonly used procedure \cite{CY,Baez,I1}.
The concept of an additive category generalizes that of a ring
in the same way groupoids generalize the notion of groups.
Additive categories were called "rings with several objects"
in \cite{Mi}, and were studied by imitating
results and proofs from noncommutative homological ring theory,
to additive category theory.
Alternatively, the additive category theory may be applied, as in
\cite{St}, to ring theory.
Subsequent related papers adopted the "ideal theory" point
of view, e.g. \cite{Ch}, and in \cite{SD} the problem of
lifting algebraic geometry to the category theory level
was considered and a notion of prime spectrum of a category
was defined.

In this article we consider the Dedekind's original aim for introducing
ideals \cite{Ed}, and leading to the study of general rings,
not only principal ideal rings (PIR).
In the context of categories, we relax the requirements
of an exact category for the existence kernels and cokernels, 
and define homological objects in an intrinsic way, using ideals.
The former are not ``intrinsic'' concepts.
Universal constructions in category theory represent classes of morphisms
in terms of a (universal) generator (a limit).
The use of ``coordinates'' in geometry, 
or the use of generators and relations in algebra is opposed to
the intrinsic point of view which emphasizes the ``global object'',
which in our case is the ideal.

To fix the ideas, consider an example, kernels:
$$\diagram
A \rto^{f} & B\\
K \ar@{^{(}->}[u]^{ker(f)} &
X \ar@{..>}[l]^{\exists !\phi} \ulto_{x} \uto_{0}
\enddiagram$$
Via categorification, the morphisms whose composition
with $f$ is zero is in the one object category case (ring)
a right ideal.
The existence of kernels in the categorical sense means
that this right ideal is principal and the universality property
means it is free as a right module, i.e. that the one generator
forms a base, and therefor it is unique modulo a ``unit'' (isomorphism):
$$Ker(f):=\{ker(f)\circ\phi | \forall \phi\}\ {\bf =}\
Ann_R(f):=\{x|fx=0\}$$
$$ker(f) \ \text{exists (categorical level) iff}
Ann_R(f)\ \text{is a free \& {\bf principal} right ideal
(ring theory level)}$$
What benefit we have from this reinterpretation?
To define homology of complexes in a category, one normally
assumes that there exists a null family of morphisms and
kernels and cokernels exist.
Instead, one may consider the intrinsic approach to generalize homology
in the context of additive categories, using ideals.

\vspace{.2in}
In section \ref{S:ideals} we review the notion of an ideal in a category,
and introduce some new related concepts as well as operations with
ideals.
In  section \ref{S:homol} we define the notions of kernel, cokernel,
image and coimage in terms of ideals.
These allow the definition of the homology of
a complex in an additive category in an intrinsic way.
It is not just a representable version of the ordinary homology,
unless the category is exact. Then the homology modules are
represented on projectives by the usual homology groups
(theorem \ref{T:rep1}).
Applications to derived categories of abelian categories which are not exact,
are considered in section \ref{S:appl}.
The construction of derived functors of functors which are not additive
is considered.
In section \ref{S:axioms} we consider axioms for an additive category,
formally yielding the structure of an abelian category
in terms of ideals.

\section{Ideals in additive categories}\label{S:ideals}
Let $\C{A}$ be a additive category, i.e. a category such that
for any two objects $A, B$ in $\C{A}$, $Hom(A,B)$
is an object in the category $\C{A}b$ of abelian groups.
We consider the notion of ideal as originaly defined by \cite{Ke} (p.300),
and consistent with the correspondence rings-additive categories,
through categorification.
\begin{defin}
An ideal (left/right) is a collection
$\{I(A,B)\}_{A,B\in Ob(\C{A})}$ of abelian subgroups
$I(A,B)\subset Hom_\C{A}(A,B)$, which is stable under
(left/right) composition  with morphisms (whenever defined).

We will consider also such families of abelian subgroups indexed by 
objects of subcategories of $\C{A}$.
If $\C{S}, \C{T}$ are subcategories of $\C{A}$,
and ideal (left/right) from $\C{S}$ to $\C{T}$ is a family 
$\{I(A,B)\}_{A\in Ob(\C{S}), B\in Ob(\C{T})}$ of abelian subgroups
$I(A,B)\subset Hom_\C{A}(A,B)$, which is stable under
composition (left/right) with morphisms from $\C{A}$.
\end{defin}
An ideal (bilateral) is a left ideal which is also a right ideal.
Ideals are in a one to one correspondence with
subfunctors of $Hom:\C{A}^{op}\boxtimes \C{A}\to \C{A}b$ (\cite{ML}, p.36).

Left (right) ideals may be viewed as families $I(A)_{A\in Ob(\C{A})}$
of subfunctors of the corresponding canonical representable functors
$Hom(A, \cdot)$.
The notation $I(A)=I_A$ will be also used.

The condition that $I(A)$ is a subfunctor of $Hom(A, \cdot)$
ensures that the family is stable under left composition with
arbitrary morphisms $\phi:X\to Y$, where $A, X, Y\in Ob(\C{A})$:
$$\diagram
I(A)(\phi)=Hom(A,\phi)_{|I(A,X)}                   &             &
\C{A} \ar@/^1pc/[r]^{I(A)} \ar@/_1pc/[r]_{I(B)} \ar@{}[r]|{\Downarrow I(f)} & \C{A}b\\
A \ar@{=>}[dr]_{I(A,Y)} \ar@{=>}[r]^{I(A,X)} & X \dto^{\phi} & I(B,X) \dto_{\phi\circ \cdot} \rto^{\cdot\circ f} & I(A,X) \dto^{\phi\circ\cdot} \\
 					     & Y             & I(B,Y) \rto_{\cdot\circ f}& I(A,Y)
\enddiagram$$
The double arrow denotes a set of morphisms, and the left diagram
is 2-commutative, i.e. $\eta_{A,X}$ is the 2-morphism corresponding
to the inclusion of sets.
Composing $\phi$ with morphisms from $I(A,X)$
may result in a proper subset of $I(A,Y)$.

The diagram from the right represents the naturality of the
transformation $I(f)$ associated to a morphisms $f:A\to B$.
Its commutativity is equivalent to the associativity of the composition
of morphisms in $\C{A}$.

There is an alternative terminology for right (left) ideals.
A right (left) ideal (from $\C{A}$) to $P\in Ob(\C{A})$ is called in
{\em right (left) ideal} in \cite{Mi} (p.18), and
a {\em right P-ideal} in \cite{St} (p.140).
It is essentially a {\em (co)sieve at $P$} \cite{ML} (p.37), by adjoining
null morphisms ending at $P$ if necesary
(see also\cite{Ba}, p.171 and \cite{SD}, p.139).

The {\em support} of the left (right) ideal $I$ is defined as the full
subcategory of $\C{S}$ consisting of objects $A$ such that
$I(A,X)\ne 0$ ($I(X,A)\ne 0$) for at least one object $X$ of $\C{A}$.
Since any ideal can be extended trivially outside its support,
we will assume all ideals defined on $\C{A}$.

As an example, for a given morphism $f$ consider the class of morphisms
$g\circ f$ ($f \circ g$) which left (right) factor through $f$.
It is the {\em principal} left (right) ideal generated by a morphism
$f$ \cite{SD} (p.142),
and it will be denoted by $<f|$ ($|f>$), or alternatively by
$\C{A}f$ ($f\C{A}$).
As another example, consider the category $\C{O}(X)$ of open sets and
inclusions, canonically associated to a topological space $X$.
The principal left ideal generated by an open set $U$
is  essentially the family of open sets containing $U$ \cite{ML} (p.70).

An ideal $I$ supported on $\C{S}$ is called {\em proper}
if it is different from $Hom_{|\C{S}}$ and $0_{|\C{S}}$
(compare \cite{SD} (p.145)).
In particular, the left ideal supported at $A$ consisting of
all morphisms with codomain the object $A$ is not a proper ideal.
It is the {\em total sieve on $A$} (maximal sieve on $A$ \cite{ML} (p.38)).
A {\em maximal ideal} is defined as usual as a proper ideal
not contained in another proper ideal \cite{SD} (p.141).

The {\em product} of two ideals (left / right) $I$ and $J$ is the
ideal generated by the class of products of morphisms $f\circ g$
with $f\in I$ and $g\in J$ (zero if the class of generators is empty).

The {\em intersection} of two left (right) ideals
is the class of common morphisms.
As an example, if the domain of $g$ is the codomain of $f$, then
the product of the left principal ideal (l.p.i.) $<g|$ and
the r.p.i. $|f>$ is a bilateral ideal equal to their {\em intersection}.
It is the bilateral ideal $<g\circ f>$ generated by their composition,
and it is denoted by $<g|f>$:
$$A\overset{f}{\to} B\overset{g}{\to} C,
\qquad <g|f>=<g|\cap |f>=<g|\cdot|f>=\C{A}(g\circ f)\C{A}$$
%
If $I$ is a bilateral ideal $\hat{\C{A}}=\C{A}/I$ is the category
having the same object as $\C{A}$ and morphisms are cosets
modulo the ideal $I$. Note that $I$ also defines a subcategory
$\C{I}$ by restricting the morphisms to those in $I$.
If $I$ and  $J\subset I$ are two ideals,
then $I/J:\C{I}/J\boxtimes\C{I}/J\to \C{A}b$ is the
ideal associating to a pair of objects $(A,B)$ the
quotient group $I(A,B)/J(A,B)$.

To define a quotient of left (right) ideals we need the notion
of a module \cite{Ba} (p.171), \cite{St} (p.140).
A left (right) $\C{A}-module$ is an additive functor
$M:\C{A}\to \C{A}b$ ($M:\C{A}^{op}\to \C{A}b$).
It extends in a natural way the one object case,
when $\C{A}$ is the one object additive category $\{*_R\}$
associated to a ring $R=End(*_R)$
(tautological categorification - \cite{I1}).
Then an $R$-module structure $\phi:R\to End(M)$ on the abelian group $M$
is the same as a functor $M:\C{A}\to \C{A}b$.
We denote by $\C{A}-mod$ ($mod-\C{A}$) the category
of left (right) $\C{A}-modules$.

%
%
\begin{defin}
If $I$ and $J\subset I$ are two right ideals, then
their {\em quotient} $I/J:\C{A}\to mod-\C{A}$ is defined as
the family of right modules $(I/J)_A=I_A/J_B$.
Denoting for clarity $F=I_A, G=J_A$ and $\eta:F\to G$ the corresponding
inclusion, we have the following diagram defining the functor $I/J$ on
morphisms ($f:X\to Y$):
$$\diagram
F(X) \dto_{F(f)} \ar@{^(->}[r]^{\eta_X} &
G(X) \dto^{G(f)} \ar@{->>}[r]^{can_X} &
F/G(X) \ar@{..>}[d]^{\exists! F/G(f)} \\
F(Y) \ar@{^(->}[r]^{\eta_Y} & G(Y) \ar@{->>}[r]^{can_Y} &
F/G(Y)
\enddiagram$$
\end{defin}
Following the non-commutative ring theory, for example \cite{Ja},
we define the {\em left annihilator} $\C{L}(S)$ of
a non-empty class of morphisms $S$ as the class of morphisms $f$
such that $f\circ a=0$ for all $a\in S$.
%
A left ideal $I$ is called an {\em annihilator left ideal} of $\C{A}$
if $I=\C{L}(S)$ for some non-empty class of morphisms $S$.

\section{Homology}\label{S:homol}
To define the homology of a complex in $\C{A}$ we consider
the families of morphisms, rather then the limits giving
the usual universal objects.
\begin{defin}
The right (left) annihilator of a left (right) ideal $I$ 
is called the {\em kernel ideal} ({\em cokernel ideal}),
and will be denoted by $Ker(I)$ ($Coker(I)$).
The {\em image ideal} is $Im(I)=Ker(Coker(I))$ and
the {\em coimage ideal} is $Coim(I)=Coker(Ker(I))$.

If $I$ is a principal left (right) ideal generated by $f$,
the following brief notation will be used:
$$Ker(f)=\C{R}(<f|), \qquad Coker(f)=\C{L}(|f>)$$
$$Im(f)=Ker(Coker(f)), \qquad Coim(f)=Coker(Ker(f))$$
\end{defin}
From the definitions the following lemma is clear.
%
%
\begin{lem}
$Ker(I)$ and $Im(I)$ are right ideals.
$Coker(I)$ and $Coim(I)$ are left ideals.

If $I=<f|$, then $Ker(I)$ is a sieve on the domain of $f$
and $Im(I)$ is a sieve on the codomain of $f$.
Also $Coker(I)$ is a cosieve on the codomain of $F$ and
$Coim(I)$ is a cosieve on the domain of $f$.
\end{lem}
If $Ker(f)$ ($Im(f)$) is a principal right ideal,
then a generator is a kernel of $f$ in the usual sense,
and will be denoted by $ker(f)$ ($im(f)$).
Similarly for the principal left ideals
$Coker(f)=|coker(f)>$ and $Coim(f)=|coim(f)>$.

In what follows we will  use the terms kernel, cokernel, etc.
as an abbreviation for the corresponding ideals, not assumed to be principal.
The following lemma is proved by inspection.
\begin{lem}
If the morphisms $g$ and $f$ are compassable, then the following conditions
are equivalent:

a) The sequence of two morphisms is a complex: $g\circ f=0$.

b) The image of $f$ is contained in the kernel of $g$:
$Im(f)\subset Ker(g)$.

c) The cokernel of $f$ contains the coimage of $g$:
$Coker(f)\supset Coim(g)$.
\end{lem}
%
%
Let $\C{C}h(\C{A})$ be the category of homological complexes in $\C{A}$.
%
\begin{defin}\label{D:homol}
The {right homology} of the complex $C_\bullet$ is
the family of right modules $\C{H}^R_\bullet(C_\bullet)$
defined as the quotient of right ideals:
$$\C{H}^R_n(C_\bullet)=Ker(d_n)/Im(d_{n+1}).$$
The {left homology} of the complex $C_\bullet$ is
the family of left modules $\C{H}^R_\bullet(C_\bullet)$
defined as the quotient of left ideals:
$$\C{H}^L_n(C_\bullet)=Coker(d_{n+1})/Coim(d_n).$$
The complex is {\em exact} iff all homology modules vanish.

Let $f_\bullet:C_\bullet\to D_\bullet$ be a morphism of complexes.
Define $\C{H}^R_n(f_\bullet): \C{H}^R_n(C_\bullet)\to \C{H}^R_n(D_\bullet)$
the natutal transformation induced by $f_\bullet$ on homology
(see remark \ref{R:homol1}).
\end{defin}
We thus obtain the functors $\C{H}^R_n:Comp(\C{A})\to Hom(\C{A}, \C{A}b)$
and $\C{H}^L_n:Comp(\C{A})\to Hom(\C{A}^{op}, \C{A}b)$ with values:
\begin{align}
H_n^R(C_\bullet)(X)=&\{\phi:X\to C_n| d_n\circ \phi\}/\notag\\
& \{\phi:X\to C_n| Coker(d_{n-1})\circ \phi=0\}\notag
\end{align}
equivalence classes of morphisms.
\begin{rem}\label{R:homol1}
It is imediate that $Ker(d_n)(X)=ker(Hom(X, d_n))$, but in general
$Im(f)(X)\ne im(Hom(X, f))$ (see example \ref{E:homol1}).
As an imediat consequence, a chain map $f_\bullet$ as above,
restricts to a map between kernels.
\end{rem}
\begin{rem}
Associate to a homological complex $C_\bullet$ the right ideal
$I(C_\bullet)$ determined by the class of morphisms consisting of
the morphisms of the complex together with the total sieves on the objects
of the complex.
Similarly, consider $J(C_\bullet)$ the left ideal obtained by
adjoining the cosieves supported on the complex.
Then, the homology may be defined ``globally''
as $Ker(J(C_\bullet))/Im(J(C_\bullet))$.
The direct definition is used for simplicity.
\end{rem}
With this notion of exactness, $Hom$ in not left exact in general.
The usual sequences stating that $Hom$ is left exact reduce to
the following sequences.
\begin{prop}\label{P:hom1}
If $0\to A\overset{f}{\to}B\overset{g}{\to}C\to 0$ is an
exact sequence in the additive category $\C{A}$,
then for any object $X\in \C{A}$, the following sequences are exact:
$$0\to Hom(X,A)\overset{f_*}{\to}Hom(X,B)\qquad
0\to Hom(C,X)\overset{g^*}{\to}Hom(B,X)$$
where $f_*=Hom(X, f)$ and $f^*=Hom(g,X)$.
\end{prop}
\begin{pf}
They are equivalent to $f_*$ and $g^*$
being monomorphisms ($Ker(f_*)=0$ and $Ker(g^*)=0$).
Indeed $f_*(\psi)=0$, i.e. $f\circ\psi=0$, is equivalent to
$\psi\in Ker(f)$, therefor $\phi=0$.
Similarly, $g^*(\phi)=0$, i.e. $\phi\circ g=0$, is equivalent to
$\phi\in Coker(g)$, therefor $\phi=0$.
\end{pf}
\begin{example}\label{E:homol1}
A typical example showing how $Hom$ may fail to be left exact,
is the following.

Let $\C{A}$ be the category of topological abelian groups and
$k$ a commutative ring, e.g. $k=\B{Z}$.
Consider the category $\C{A}=k\C{C}$ with the same objects as $\C{C}$ and
with morphisms the free $k$-modules generated by the morphisms of the
category $\C{C}$.
Then $\C{A}$ is an additive category in an obvious way.
The sequence $0\to \B{Q}\overset{j}{\to} \B{R}\to 0\to 0$
is exact in $\C{A}$, since $Coker(j)=0_\B{R}$,
so that $Im(j)=Ker(0_\B{R})=1_{\B{R}}\C{A}$.
Applying $Hom(\B{Z},\cdot)$ one gets essentially the same sequence
but in $\C{A}b$, where it is not exact.
\end{example}
%
%
The relation to homology groups when the category is abelian is
investigated next.
\begin{th}\label{T:rep1}
Let $C_\bullet$ be a complex in the abelian category $\C{A}$.
Then the right homology modules $\C{H}^R_n(C_\bullet)$ are represented
by the corresponding homology groups
on the class of projective objects of $\C{A}$.
\end{th}
\begin{pf}
The category being abelian, the kernel ideal of $d_n$ and the
image ideal of $d_{n+1}$ are principal and generated by the morphisms
shown in the following commutative diagram:
$$\diagram
\cdots C_{n+1} \rrto^{d_{n+1}} \ar@{->>}[dr]_{} & & C_n \rrto^{d_n} & & C_{n-1} \cdots\\
 & B_n \ar@{^{(}->}[ur]^{i_n} \ar@{^{(}->}[rr]^{h_n} & & Z_n \urto_{0} \ar@{^{(}->}[ul]_{j_n}\\
 & & A \ulto^{\psi} \urto_{\phi} &
\enddiagram$$
If $A$ is an arbitrary object, then  $Ker(d_n)(A)$ consists of the
morphisms $j_n\circ\phi$ with $\phi$ arbitrary in $Hom(A,Z_n)$,
where $j_n=ker(f)$. Similarly $Im(d_{n+1})(A)$ is parametrized by
$\psi\in Hom(A,B_n)$, where $i_n=im(d_{n+1})$.
Applying $Hom(P,\cdot)$ to the short exact sequence
$0\to B_n\to Z_n\to H_n(C_\bullet)\to 0$,
with $P$ projective, one obtains the following s.e.s.:
$$0\to Hom(P,B_n)\to Hom(P,Z_n)\to Hom(P, H_n(C_\bullet))\to 0.$$
Note that $\C{H}^R_n(C_\bullet)(P)=$ is isomorphic to
$Hom(P,Z_n)/Hom(P,B_n)$.
In this way $\C{H}^R_n(C_\bullet)$ is isomorphic to
$Hom(\cdot, H_n(C_\bullet))$ and $H_n(C_\bullet)$ represents
the right homology module $\C{H}^R_n(C_\bullet)$,
on the class of projective objects of $\C{A}$.
\end{pf}
%
%
In particular when the category has a projective generator $U$,
then the category is canonically embedded in the category of abelian
groups, and the two definitions of homology functors canonically correspond.

In the following diagram $Y$ denotes the Yoneda embeding and $h_U$ is
the canonical representable functor which embeds the additive category
$\C{A}$ into the abelian category of abelian groups.
Recall that a projective generator can be characterized by
the additive functor $h_U$ being exact and faithfull.
$$\diagram
{\large Comp(\C{A})} \ar@{..>}[d]_{H_n}^{\quad \overset{\eta}{\swarrow}\ \ }
\drto^{H_n^R}\\
\C{A} \ar@{_(->}[dr]_{h_U} \ar@{^(->}[r]_{Y} &
{\large Hom(\C{A}, \C{A}b)} \dto^{<\cdot, U>}\\
& \C{A}b
\enddiagram$$
As an example, consider $\C{A}=R-mod$ ($R$ commutative ring).
\begin{cor}
If $U=R$ is the canonical projective generator of $R-mod$ then
there are canonical isomorphisms:
$$H_n^R(\cdot)(U)\cong h_U\circ H_n\cong H_n$$
\end{cor}
This is essentially due to the identification of the elements of
an $R$-module $M$ with the morphisms $Hom(R,M)$.

\section{Applications}\label{S:appl}
Important examples are derived categories of abelian categories,
which are not necessarily abelian.
One would still want to have the notions of kernel and cokernel, etc.
in a generalized sense, and to study these categories using
``abelian techniques'', e.g. the machinary of derived functors.

%
%
\subsection{Derived categories}
Recall that, given an abelian category $\C{A}$, one can consider the
category of complexes and chain maps $Comp(\C{A})$ which is also abelian.

Define the category $K(\C{A})=Comp(\C{A})/\C{I}$ with the same objects
as $Comp(\C{A})$ (complexes) and with morphisms {\bf homotopy classes} of
chain maps, i.e. classes of morphisms modulo the ideal
$\C{I}\subset Hom_{Comp(\C{A})}$ of null chain homotopic chain maps.
$$\diagram
\C{A}\rrto^{C_0\quad} \drrto_{1_\C{A}} & &
{\large Comp(\C{A})} \ar@{->>}[rr] \dto^{H_0} & &
K(\C{A}) \dllto^{\hat{H}_0} \\
& & \C{A} &  &
\enddiagram$$
In the above diagram $C_0$ embeds the objects and morphisms as complexes
and chain maps concentrated in degree 0.

We will not need the actual derived category
$\C{D}(\C{A})=K(\C{A})[\Sigma^{-1}]$, which is obtained as a localization
with respect to quasi-isomorphisms.

In the general case, the category $K(\C{A})$,
obtained by considering homotopy classes
of morphisms, is no longer abelian.
It may fail even to be exact, as shown by the next example
\cite{Bo}, p.45.

If $\C{A}$ is a category with nontrivial extensions
(e.g. abelian groups), and
if $u\in Hom_\C{A}(X,Y)$ is a morphism with its image $Im(u)$ not a
direct summand in $Y$, then the corresponding morphisms in $Comp(\C{A})$
(concentrated in degree 0) has no ``classical'' cokernel:
$$\diagram
X_\bullet \rto^{u_\bullet} &
Y_\bullet \ar@{_(->}[dr]_{q} \rto_?^{coker(u_\bullet)}&
Z_\bullet \ar@{..>}[d]_{\exists}^{\uparrow u'}\\
 & & Cone(u) &
\enddiagram$$
The key is that $q\circ u$ is chain homotopic to zero, and
if $u$ would have a cokernel, then the s.e.s. corresponding
to its image (in $\C{A}$) must be split.

Note that the {\em ideal} $Coker(u)$ always exists.
%
%
\subsection{Derived functors}
Since the key lemmas ``lift'' from generators to ideals, 
homology theories can be defined for additve categories.


There is an analog of the connecting transformation, which in the
abelian case is represented by the usual connecting morphism.
\begin{lem}
Restricting to projective objects and morphisms, there are natural
transformations $\delta_n^R$ such that $(H_n^R, \delta_n^R)$ is a
connected sequence of functors.
\end{lem}
Similarly, the usual lemma leading to the machinery of
derived functors holds.
\begin{lem}Chain homotopic morphisms induce canonical isomorphims in
homology:
$$\quad f_\bullet\sim g_\bullet \Rightarrow
\C{H}^R_n(f_\bullet)=\C{H}^R_n(g_\bullet).$$
\end{lem}
%
%
This approach allows one to generalize the construction of derived functors
to the case of nonabelian functors.

More precisely, for any category one may consider the category with
the same objects, and with $Hom(X,Y)$ the free module generated by
the corresponding morphisms of the original category.
Enriching the category in this way (``linearization'') may be viewed as
the categorification of the group ring construction.
In this way one may apply the machinery of derived functors
in the general case of an arbitrary category.

\section{Additive categories and abelian category axioms}\label{S:axioms}
We will briefly consider the abelian category axioms \cite{St} (p.83)
in terms of ideals.
Recall that we use the term kernel (etc.) as an abbreviation
for the kernel ideal.

The analog of the axioms (K) and (K$^{op}$) of an abelian category,
hold by definition.
\begin{lem}
Let $\C{A}$ be an additive category. Then:

1) ('K). Every morphism $f$ has a kernel: $Ker(f)$.

2) ('K$^\circ$). Every morphism $f$ has a cokernel: $Coker(f)$.

3) $f$ is a monomorphism iff $Ker(f)=0$.

4) $f$ is an epimorphisms iff $Coker(f)=0$.
\end{lem}
\begin{pf}
The statements are immediate consequences of the definitions.
\end{pf}
Note that a left (right) ideal $I$ is contained in its coimage (image).
A left (right) ideal $I$ will be called {\em closed} if
it coincides with its coimage (image).
\begin{lem}
Kernels, cokernels, images and coimages are closed ideals.
\end{lem}
\begin{pf}
We will prove the statement corresponding to kernels. Note first
that taking kernels reverses the order relation defined by inclusion.
If $I$ is a left ideal, then $I\subset Coker(Ker(I))$, so that
$Ker(Coker(Ker(I)))\subset Ker(I)$. Now, if $\psi\in Ker(I)$,
for all $\psi\in Coker(ker(I))$\ $\psi\circ\phi=0$, and
$\phi\in Ker(Coker(Ker(I)))$. Therefor $Ker(Coker(Ker(I)))=Ker(I)$.
\end{pf}
The analog of the axioms (N) and (N$^\circ$) are:

\ ('N) For every monomorphism $f$, the right ideal $|f>$ is
a kernel ideal ($|f>=Ker(I)$ for some left ideal $I$).

\ ('N$^\circ$) For every epimorphism $f$, the left ideal $<f|$ is
a cokernel ideal ($<f|=Coker(I)$ for some right ideal $I$).

Note that in an abelian category the above statements reduce
to the usual axioms.
\begin{lem}
If $\C{A}$ is an abelian category then (N) and (`N) are equivalent.
Dually, (N$^{op}$) and ('N$^{op}$) are equivalent.
\end{lem}
\begin{pf}
If $f\C{A}$ is a kernel ideal $Ker(I)$, and therefor closed,
$f\C{A}=Ker(J)$ with $J=Coker(Ker(I))$ the closure of $I$.
Then $J=Coker(f\C{A})$ is a principal left ideal  generated by
$h=coker(f)$ and $f\C{A}=Ker(\C{A}h)=Ker(h)=ker(h)\C{A}$. Now it is
easy to see that two monomorphisms, $f$ and $g=ker(h)$,
generate the same right ideal iff they have isomorphic domains
$f=g\circ s$ ($s$ isomorphism).
Therefor $f$ is the kernel (in the usual sense) of a morphism.
\end{pf}
\begin{rem}
In terms used in ring theory \cite{Ja} (p.2),
the axiom ('N) (dually ('N$^{op}$))
can be rephrased as ``any {\em right regular element} $f$
(no right zero-divisors / monomorphism)
generates an {\em annihilator right ideal} $f\C{A}$''.
It further demonstrates the similarity between ring theory and
additive category theory, as the ``ring with several objects'' case.
\end{rem}
\begin{prop}
('N) holds iff any monomorphism $f$ generates its image ideal: $|f>=Im(f)$.

('N$^\circ$) holds iff any epimorphism $g$ generates its
coimage ideal: $<g|=Coim(g)$.
\end{prop}
\begin{pf}
We will only prove the first statement,
since the second statement follows by duality.
If $|f>=f\C{A}$ is a kernel ideal $Ker(I)$, then it is closed and
$Im(f)=Ker(Coker(f\C{A}))$ equals $f\C{A}$.

For the converse, let $f$ be a monomorphism, and assume $f\C{A}=Im(f)$.
Then $f\C{A}=Ker(I)$ with $I=Coker(f)$.
\end{pf}
It is not clear what the substitute for epi-mono factorization axiom
(AB2) \cite{St} (p.84) should be,
and allowing the usual epi-mono based type of arguments.
We only note that the rows in the next diagram are exact,
i.e. the kernel of the left ideal based at $A$ equals
the image of the right ideal based at $A$:
$$\diagram
\ar@{=>}[r]^{Ker(f)} & A \dto_{f} \ar@{=>}[r]^{Coim(f)} & & Ker(Coim(f))=Im(Ker(f)\\
\ar@{<=}[r]^{Coker(f)} & B \ar@{<=}[r]^{Im(f)} & & Ker(Coker(f))=Im(Im(f)\\
\enddiagram$$
and that in general $f$ is common to both $Im(f)$ and $Coim(f)$.

\section{Conclusions and further developments}
The aim of ``categorifying ideal theory'' is manifold.
First, to avoid the familiar "unique-up-to-a-unit" of arithmetic in
principal ideal rings and to define an "intrinsic" homology
in additive categories.
Second, with ideals replacing generators, a greater degree of generality
is achieved.
In additive categories, the analog of two basic axioms in
abelian categories are consequences of the definitions.

In this way it is possible to develop the usual formalism of
homological algebra in derived categories, based on short exact sequences.
Then a natural problem is to determine the relationship between this
approach and the established formalism of distinguished triangles.




\begin{thebibliography}{33}

\bibitem[Baez]{Baez} Baez, J. C., Dolan, J.,
Categorification,
{\em preprint math.QA/9802029}, (1998).

\bibitem[Ba]{Ba}
Baumgartner, Karlbeinz,
Structure of additive categories,
{\em Cahier de Topologie et geometrie differentielle}
Vol. {\bf XVI-2} (1975), pp.169-213.

\bibitem[Bo]{Bo}
Borel, A.  et. al.,
Algebraic D-modules,
Perspectives in Mathematics, Vol. {\bf 2},
Academic Press, Inc., 1987.

\bibitem[Ch]{Ch}
Christensen, Daniel, J.,
Ideals in triangulated categories: phantoms, ghosts and skeleta,
{\em math.AT/9807071} (1998).

\bibitem[CY]{CY} Crane, L., Yetter, D.
Examples of categorification,
{\em preprint math.QA/9607028}, (1996).

\bibitem[Ed]{Ed}
Edwards, Harold M.,
Mathematical Ideas, Ideals, and Ideology,
{\em The Mathematical Inteligencer} Vol. {\bf 14}, No. 2, 1992,
pp. 6-19.

\bibitem[I1]{I1}
Ionescu, Lucian M.,
On categorification,
math.CT/9906038.

\bibitem[Ja]{Ja}
Jategaonkar, A. V.,
Left Principal Ideal Rings,
{\em Lecture Notes in Mathematics} {\bf 123},
Springer-Verlag, Berlin Heidelberg New York 1970.

\bibitem[Ke]{Ke}
Kelly, G., M.,
On the radical of a category,
{\em J. Australian Math. Soc.} {\bf 4} (1964) pp.299-307.

\bibitem[Mi]{Mi}
Mitchell, B.,
Rings with several objects,
Advances in Mathematics {\bf 8} (1972), pp.1-161.

\bibitem[ML]{ML}
MacLane, Saunders,
Sheaves in Geometry and Logic,
{\em Springer-Verlag New York}, 1992.

\bibitem[SD]{SD}
Shu-Hao Sun, Ajneet Dhillon,
Prime spectra of additive categories (I),
{\em J. Pure Appl. Algebra} {\bf 122} (1997), pp.135-157.

\bibitem[St]{St}
Street, Ross,
Ideals, Radicals, and Structure of Additive Categories,
{\em Applied Categorical Structures} {\bf 3} (1995), pp.139-149.

\bibitem[Str]{Str}
Strooker, Jan R.,
Introduction to categories, homological algebra and sheaf cohomology,
{\em Cambridge University Press},
Cambridge London New York Melbourne, 1978.

\end{thebibliography}
\end{document}